\documentclass[12pt]{article}
\usepackage{mathrsfs,graphicx,amsmath,amssymb,amsfonts}
\begin{document}
\title{The effect of a graft transformation on distance signless Laplacian spectral radius of the graphs \thanks{This work is supported by NSFC (No. 11461071).}}
\author{Dandan Fan$^{a}$, Guoping Wang$^{b}$\footnote{Corresponding author. Email: xj.wgp@163.com.}, Yinfeng Zhu$^{b}$\\
{\small $^{a}$ College of Mathematical and Physical Sciences, Xinjiang Agricultural University,}\\
{\small Urumqi 830052, Xinjiang, P. R. China.}\\
{\small $^{b}$ School of Mathematical Sciences, Xinjiang Normal University,}\\
{\small Urumqi 830054, Xinjiang, P. R. China.}}

\date{}
\maketitle {\bf Abstract.} Suppose that the vertex set of a connected graph $G$ is $V(G)=\{v_1,\cdots,v_n\}$.
Then we denote by $Tr_{G}(v_i)$ the sum of distances between $v_i$ and all other vertices of $G$.
Let $Tr(G)$ be the $n\times n$ diagonal matrix with its $(i,i)$-entry equal to $Tr_{G}(v_{i})$
and $D(G)$ be the distance matrix of $G$. Then $Q_{D}(G)=Tr(G)+D(G)$ is the distance signless Laplacian matrix of $G$.
The largest eigenvalues of $Q_D(G)$ is called distance signless Laplacian spectral radius of $G$.
In this paper we give some graft transformations on distance signless Laplacian spectral radius of the graphs
and use them to characterize the graphs with the minimum and maximal distance signless Laplacian spectral radius
among non-starlike and non-caterpillar trees.

\section{Introduction}
\qquad Let $G$ be a simple connected graph with the vertex set $V(G)=\{v_1,\ldots,v_n\}$.
The {\it distance} $d_{G}(u,v)$ between the vertices $u$ and $v$ is the length of shortest path between $u$ and $v$ in $G$.
For $u\in V(G)$, the {\it transmission} $Tr_{G}(u)$ of $u$ is the sum of distances between $u$ and all other vertices of $G$.
Let $Tr(G)$ be the $n\times n$ diagnonal matrix with its $(i,i)$-entry equal to $Tr_{G}(v_{i})$ and $D(G)$ be the distance matrix of $G$.
Then $Q_{D}(G)=Tr(G)+D(G)$ is the distance signless Laplacian matrix of $G$.
The largest eigenvalues $\rho_{Q}(G)$ of $Q_{D}(G)$ is the {\it distance signless Laplacian spectral radius} of $G$.

The distance spectral radius of a connected graph has been studied extensively.
S. Bose, M. Nath and S. Paul \cite{S. Bose} determined the unique graph with maximal distance
spectral radius in the class of graphs without a pendant vertex.
G. Yu, Y. Wu, Y. Zhang and J. Shu \cite{G. Yu-2} obtained respectively the extremal graph and unicyclic graph
with the maximum and minimum distance spectral radius.
W. Ning, L. Ouyang and M. Lu \cite{W. Ning} characterized the graph with minimum distance spectral radius
among trees with given number of pendant vertices.
G. Yu, H. Jia, H. Zhang and J. Shu \cite{G. Yu} characterized the unique graphs
with the minimum and maximum distance spectral radius among trees and graphs with given number of pendant vertices.
For more about the distance spectra of graphs see the surveys \cite{M. Aouchiche-2, D. Stevanovic} as well as the references therein.

M. Aouchiche and P. Hansen introduced in \cite{M. Aouchiche} the distance Laplacian and distance signless Laplacian spectra of graphs,
and proved in \cite{M. Aouchiche-1} that the star $S_n$ attains minimum distance
Laplacian spectral radius among all trees of order $n$.
R. Xing and B. Zhou \cite{R. Xing} gave the unique graph with minimum distance and
distance signless Laplacian spectral radii among bicyclic graphs.
R. Xing, B. Zhou and J. Li \cite{R. Xing-2} determined the graphs with minimum
distance signless Laplacian spectral radius among the trees, unicyclic graphs,
bipartite graphs, the connected graphs with fixed pendant vertices and fixed connectivity, respectively.
A. Niu, D. Fan and G. Wang \cite{A. Niu} determined the extremal graph with the minimum distance Laplacian spectral radius
among bipartite graphs with the fixed connectivity and matching numbers.
H. Lin and B. Zhou \cite{H. Lin-2} characterized the unique graphs with the minimum distance Laplacian spectral radius with fixed number of pendant
vertices and edge connectivity among graphs, the unique tree with the minimum distance Laplacian spectral radius among trees with fixed bipartition.

A tree $T$ is {\it starlike} if it contains at most one vertex of degree at least $3$, and otherwise it is {\it non-starlike}.
A tree is {\it caterpillar} if the deletion of all pendant vertices yields a path, and otherwise it is {\it non-caterpillar}. Set
$$\mathfrak{T}(n)=\{T|~T~ \mbox{is a non-caterpillar tree of order~n.}\}$$
$$\Re(n)=\{T|~T~ \mbox{is a non-starlike tree of order~n.}\}$$

R. Xing, B. Zhou and F. Dong \cite{R. Xing-3} determined the graphs with the minimum and maximal distance spectral radius in $\mathfrak{T}(n)$ and $\Re(n)$ respectively.
In this paper we give some graft transformations on distance signless Laplacian spectral radius of the graphs
and then use them to characterize the graphs with the minimum and maximal distance signless Laplacian spectral radius in $\mathfrak{T}(n)$ and $\Re(n)$ respectively.

\section{Trees with minimum distance Signless Laplacian spectral radius in $\Re(n)$ and $\mathfrak{T}(n)$}
\qquad Suppose that $G$ is a connected graph with $V(G)=\{v_{1},v_{2},\ldots,v_{n}\}$.
Let $x=(x_{v_1},x_{v_2},\ldots,x_{v_n})^T$, where $x_{v_i}$ corresponds to $v_{i}$, i.e. $x(v_i)=x_{v_i}$ for $i=1, 2,\ldots, n$.
Then $$x^{T}Q_{D}(G)x=\sum\limits_{\{u,v\}\subseteq V(G)}d_G(u,v)(x(u)+x(v))^{2}$$
which shows that $Q_{D}(G)$ is positive definite for $n\geq 3$.\vskip 2mm

Suppose that $x$ is an eigenvector of $Q_D(G)$ corresponding to the eigenvalue $\lambda$. Then for each $v\in V(G)$,
$\lambda x(v)=\sum\limits_{u\in V(G)}d_G(u,v)(x(u)+x(v)).$\vskip 2mm

By the Perron-Frobenius theorem, we know that there is the unique unit vector $x>0$ satisfies $Q_{D}(G)x=\rho_{Q}(G)x$, and we call $x$ the {\it Perron vector} of $G$ corresponding to $\rho_{Q}(G)$.\vskip 2mm

Suppose that $G$ is a graph. If a vertex of $G$ is of degree one then it is a {\it pendant vertex}.
An edge and a path are respectively {\it pendant edge} and {\it pendant path} if one of their end vertices are pendant vertices.

Suppose that $uv$ is a cut edge of connected graph $G$.
Then we denote by $G_{uv}$ the graph obtained from $G$ by contracting $uv$ to a vertex $u$ and
attaching a pendant edge at $u$, and call such a process {\it C-transformation} of $G$ for the edge $uv$.\vskip 2mm

{\bf Lemma 2.1.} {\it Let $uv$ be a cut edge of connected graph $G$.
If $uv$ is not a pendant edge and there exists a pendant edge $uz$, then $\rho_{Q}(G)>\rho_{Q}(G_{uv})$.}

{\bf Proof.}  Let $x$ be the Perron vector of $G_{uv}$ corresponding to $\rho_{Q}(G_{uv})$.
By symmetry, we have $x_{v}=x_{z}$.

Suppose that $G\backslash \{uv\}=G_1\cup G_2$ where $u\in G_1$ and $v\in G_2$.

It can be easily observed that $d_G(u,w)=d_{G_{uv}}(u,w)$ for $w\in V(G_1)$ and that $d_{G}(w,v_{j})=d_{G_{uv}}(w,v_{j})+1$ for $w\in V(G_{1})$ and $v_j\in V(G_{2})\backslash \{v\}$.
Thus,
\begin{equation}
\begin{split}
&\rho_{Q}(G)-\rho_{Q}(G_{uv})\geq x^{T}(Q_{D}(G)-Q_D(G_{uv}))x\\
&~~~~~~~~~~~~~~~~~~~~~~~\geq\sum\limits_{w\in V(G_{1})}\sum\limits_{v_{j}\in V(G_{2})\backslash \{v\}}(x_{w}+x_{v_j})^{2}
\!-\!\sum\limits_{v_{j}\in V(G_{2})\backslash \{v\}}(x_{v}+x_{v_j})^{2}\\
&~~~~~~~~~~~~~~~~~~~~~~~=\!\sum\limits_{w\in V(G_{1})\backslash \{z\}}\!\sum\limits_{v_{j}\in V(G_{2})\backslash \{v\}}(x_{w}+x_{v_j})^{2}\!+\!\sum\limits_{v_{j}\in V(G_{2})\backslash \{v\}}(x_{z}+x_{v_j})^{2}\\
&~~~~~~~~~~~~~~~~~~~~~~~-\sum\limits_{v_{j}\in V(G_{2})\backslash \{v\}}(x_{v}+x_{v_j})^{2}>0.
\nonumber
\end{split}
\end{equation}

Suppose on the contrary that $\rho_{Q}(G_{uv})=\rho_{Q}(G)$. Then $x^{T}Q_D(G)x=\rho_{Q}(G_{uv})$.

Therefore, $(\rho_{Q}(G)-\rho_{Q}(G_{uv}))x_{u}=\sum\limits_{v_{j}\in V(G_{2})\backslash \{v\}}(x_{u}+x_{v_j})>0$.
This contradiction shows that $\rho_{Q}(G)>\rho_{Q}(G_{uv})$.
$\Box$

{\bf Lemma 2.2.} {\it Suppose that $n=n_{0}+n_{1}+n_{2}+n_{3}+4$ with $max\{n_i|i=1, 2, 3\}>1$.
If $G$ and $G'$ are shown in Fig. $1$, then $\rho_{Q}(G)>\rho_{Q}(G')$.}
\begin{center}
\begin{figure}[htbp]
\centering
\includegraphics[height=30mm]{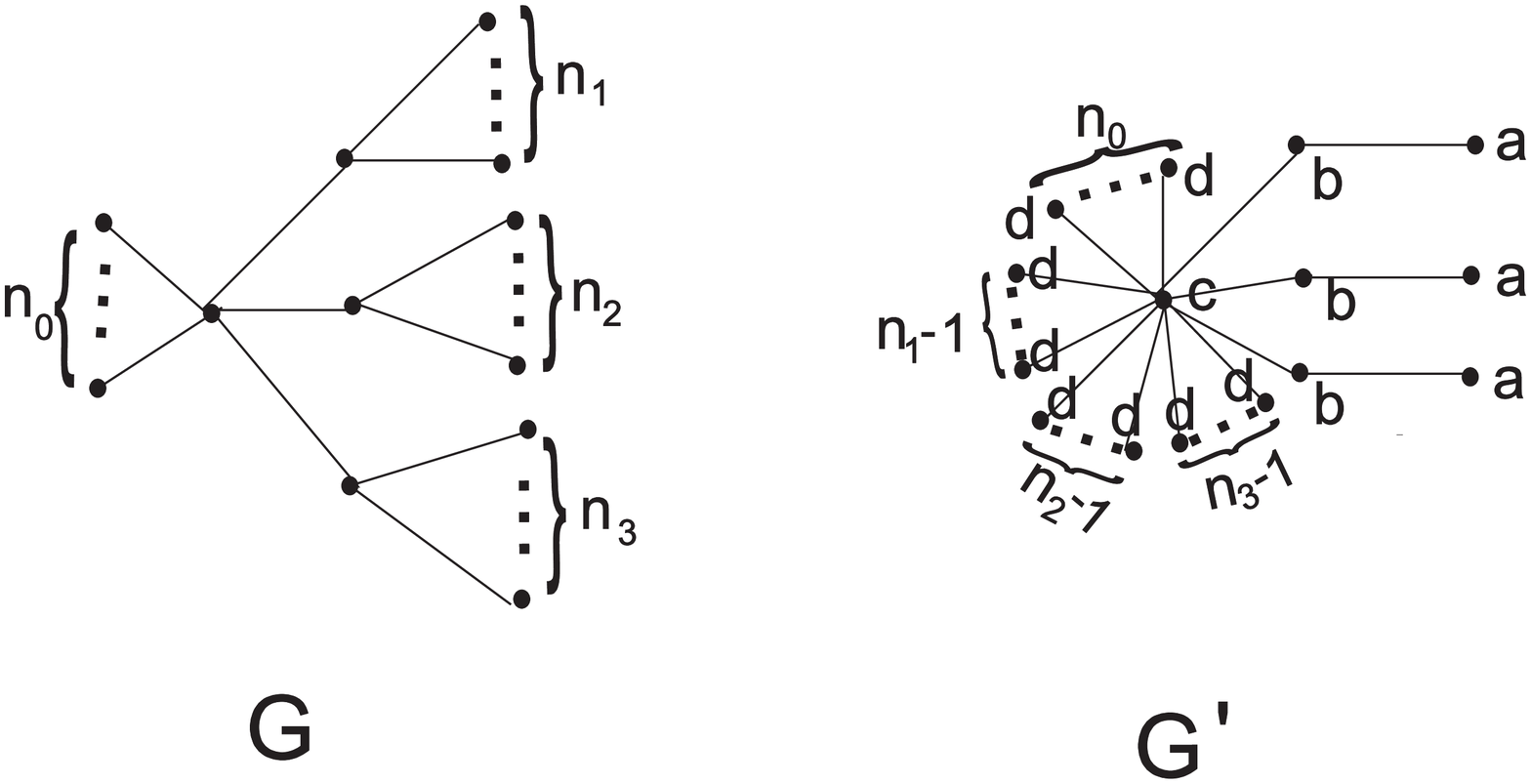}
\\ Fig. 1. ~~$G$~and~$G'$
\end{figure}
\end{center}
{\bf Proof.} Let $x$ be the Perron vector of $G'$ corresponding to $\rho_{Q}(G')$.
By symmetry, we can set $x=(a,a,a,b,b,b,\underbrace{d,\ldots,d}_{n-7},e)^T$ labeled in Fig. 1.
Without lost of generality, we assume that $n_{3}\!\geq\!n_{2}\!\geq\!n_{1}$.
Then $n_{3}\geq 2$, and hence

\begin{equation}
\begin{split}
&~~\rho_{Q}(G)-\rho_{Q}(G')\geq x^{T}(Q_{D}(G)-Q_{D}(G'))x\\
&~~~~~~~~~~~~~~~~~~~~~~~=\sum\limits_{i=1}^{3}(n_{i}-1)((d+a)^{2}\!+\!(d+b)^{2}\!+\!(d+c)^{2}\!+\!(2n-2n_{i}-n_{0}-12)(d+d)^{2})\\
&~~~~~~~~~~~~~~~~~~~~~~~>0. ~~~~~~~\Box\\
\nonumber
\end{split}
\end{equation}

Let $Tr_{max}(G)$ be the maximum transmission of vertices of $G$. Then we have the following\vskip 2mm

{\bf Lemma 2.3} \cite{H. Minc}. {\it Let $G$ be a connected graph. Then $\rho_{Q}(G)>Tr_{max}(G)$.}\vskip 2mm

{\bf Lemma 2.4.} {\it Suppose that $n=n_{0}+n_{1}+n_{2}+n_{3}+4$ with $n_{1}+n_{2}+n_{3}>4$.
If $G$ and $G^{*}$ the trees as in Fig. $2$, then $\rho_{Q}(G)>\rho_{Q}(G^{*})$.}
\begin{center}
\begin{figure}[htbp]
\centering
\includegraphics[height=30mm]{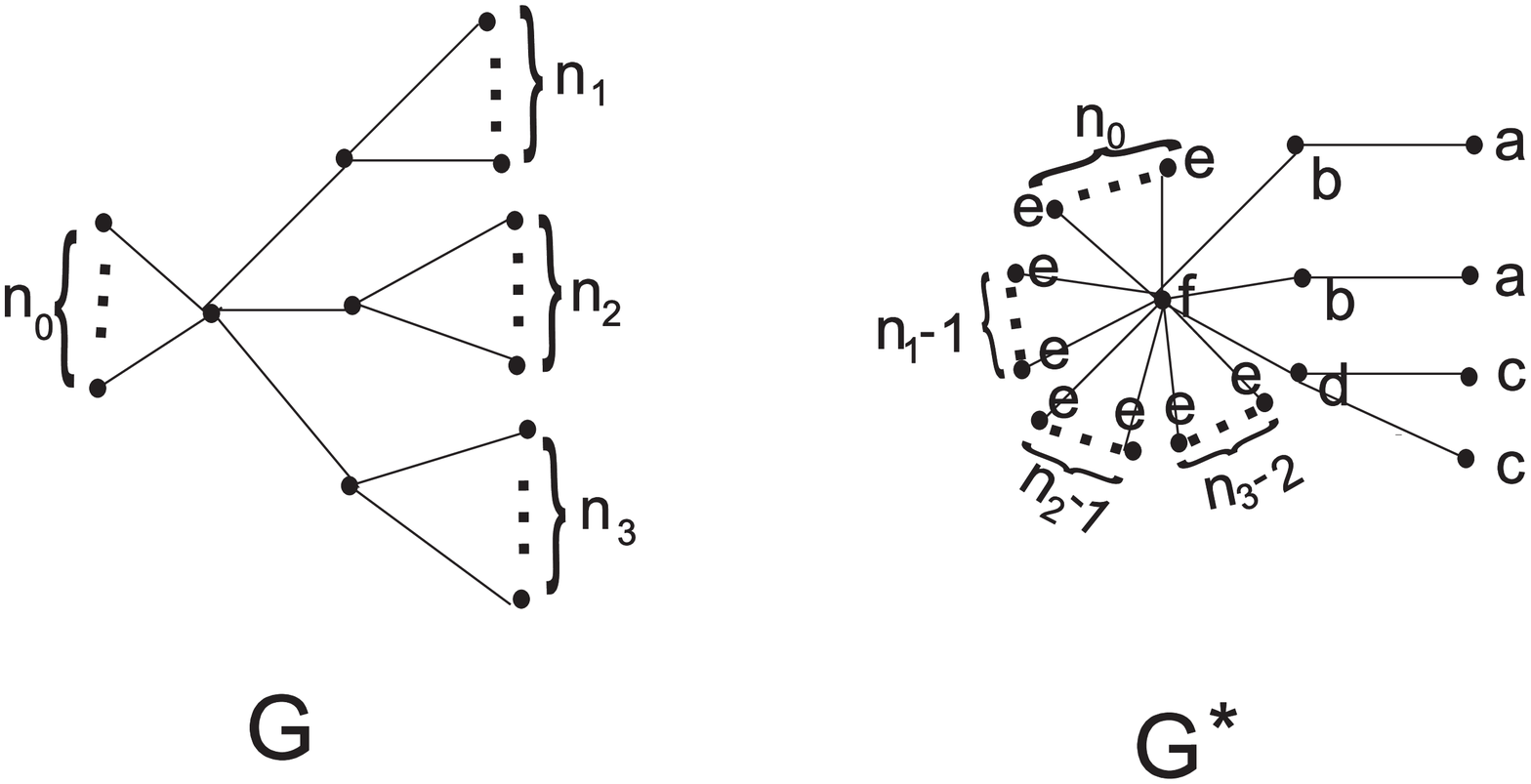}
\\ Fig. 2. ~~$G$~and~$G^*$
\end{figure}
\end{center}

{\bf Proof.} Let $x$ be the Perron vector of $G^{*}$ corresponding to $\rho_{Q}(G^{*})$.
By symmetry, we can set $x=(a,a,b,b,c,c,d,f,\underbrace{e,\ldots,e}_{n-8})^T$ labeled in Fig. 2.
Without lost of generality we can assume that $n_{3}\!\geq\!n_{2}\!\geq\!n_{1}$.

From $Q_{D}(G^{*})x=\rho_{Q}(G^{*})x$, we have

%
$$\left\{
\begin{array}{ll}
(\rho_{Q}(G^{*})-(2n-3))(a-b)=na+2c+d+(n-8)e+f,\\
(\rho_{Q}(G^{*})-(2n-4))(c-d)=2a+2b+f+(n-1)c+(n-8)e,\\
(\rho_{Q}(G^{*})-(3n-11))(a-c)=2(a-b)+2a+2d,\\
(\rho_{Q}(G^{*})-(3n-15)(2c-a)=8b+16c-d.
& \hbox{}
\end{array}
\right.$$

By Lemma 2.3, $\rho_{Q}(G^{*})>Tr_{max}(G^{*})=3n-3$ and $n>8$, and so we easily observe from the above equations that $a>b$, $c>d$, $a>c$ and $2c>a$.

We also have $(\rho_{Q}(G^{*})-(2n-5))(b-d)=2b+4c-2a$, from which we obtain $b>d$.

Thus we know that $2(e\!+\!a)^{2}\!+\!2(e\!+\!b)^{2}\!+\!(e\!+\!f)^{2}\!+\!(n\!-\!n_{3}\!-\!6)(e\!+\!e)^{2}\!-\!2(e\!+\!c)^{2}\!-\!(e\!+\!d)^{2}>0.$

It follows from $n_{1}+n_{2}+n_{3}>4$ that either $n_{1}\geq 2$ and $n_{2}\geq 2$ or $n_{3}\geq 3$.
Therefore,
\begin{equation}
\begin{split}
&\rho_{Q}(G)-\rho_{Q}(G^{*})\geq x^{T}(Q_{D}(G)-Q_{D}(G^{*}))x\\
&~~~~~~~~~~~~~~~~~~~~=\sum\limits_{i=1}^{2}(n_{i}\!-\!1)(2(e\!+\!c)^{2}\!+\!(e\!+\!d)^{2}\!+\!(e\!+\!f)^{2}
\!+\!(n\!-\!n_{i}\!-\!7)(e\!+\!e)^{2})\!+\!\\
&~~~~~~~~~~~~~~~~~~~~(n_{3}\!-\!2)(2(e\!+\!a)^{2}\!+\!2(e\!+\!b)^{2}\!+\!(e\!+\!f)^{2}\!+\!(n\!-\!n_{3}\!-\!6)(e\!+\!e)^{2}\!-\!2(e\!+\!c)^{2}\!-\!(e\!+\!d)^{2})\\
&~~~~~~~~~~~~~~~~~~~~>0.~~~~~\Box\\
\nonumber
\end{split}
\end{equation}

Let $B(n;n_{0},n_{1},\ldots,n_{r})$ be the graph obtained from the star $S_{1,r}$ by attaching $n_{0}$ edges to its center and $n_{i}$ edges to the $i$-th
pendant vertex for each $i$, where $n_{0},n_{1},\ldots,n_{r}$ are integers
satisfying $n_{0}\geq 0$, $n_{r}\geq n_{r-1}\geq \ldots \geq n_{1}\geq 1$ and $\sum\limits_{i=0}^{r}n_{i}=n-r-1$.\vskip 2mm

{\bf Theorem 2.5.} {\it If $T\in\mathfrak{T}(n)$ with $n\geq 7$, then $\rho_{Q}(T)\geq\rho_{Q}(B(n;n-7,1,1,1))$, where the equality holds if and only if $T\cong B(n;n-7,1,1,1)$.}

{\bf Proof.} Suppose that $T'$ is a tree with minimal distance signless Laplacian spectral radius in $\mathfrak{T}(n)$.
Let $d$ be the diameter of $T'$.
Now we verify that the following three claims are true.

{\bf Claim 1.} $d=4$.

Since $T'\in\mathfrak{T}(n)$, neither $d=2$ nor $d=3$.
Suppose by a way of contradiction that $d\geq 5$.
Let $P=v_{1}v_{2}\ldots v_{d+1}$ be a path of $T'$.
Making a C-transformation of $T'$ for the edge $v_{2}v_{3}$ and $v_{d-1}v_{d}$
we obtain the resulting graph $T'_{v_{2}v_{3}}$ and $T'_{v_{d-1}v_{d}}$, respectively.
Obviously, at least one of $T'_{v_{2}v_{3}}$ and $T'_{v_{d-1}v_{d}}$ falls into $\mathfrak{T}(n)$.
By Lemma 2.1, $\rho_{Q}(T')>\rho_{Q}(T'_{v_{2}v_{3}})$ and $\rho_{Q}(T')>\rho_{Q}(T'_{v_{d-1}v_{d}})$.
This contradicts minimality of $T'$, and so $d=4$.

The Claim 1 shows that $T'\cong B(n;n_{0},n_{1},n_{2},\ldots,n_{r})$, where $r\geq 3$.

{\bf Claim 2.} $r=3$.

Suppose on the contrary that $r\geq 4$. Let $T''=B(n;n_{0}+n_{r}+1,n_{1},n_{2},\ldots,n_{r-1})$.
Then $T''\in\mathfrak{T}(n)$, and by Lemma 2.1, $\rho_{Q}(T')>\rho_{Q}(T'')$.
This contradiction shows $r=3$, and so $T'\cong B(n;n_{0},n_{1},n_{2},n_{3})$.

{\bf Claim 3.} $T'\cong B(n;n-7,1,1,1)$.

If some one of $n_1$, $n_2$ and $n_3$ is greater than $1$,
then, by Lemma 2.2, $\rho_{Q}(B(n;n_{0},n_{1},n_{2},n_{3}))>\rho_{Q}(B(n;n-7,1,1,1))$,
and so $T'\cong B(n;n-7,1,1,1)$. $\Box$\vskip 2mm

{\bf Theorem 2.6.} {\it If $T \in\mathfrak{T}(n)\cap\Re(n)$ and $n\geq 8$, then $\rho_{Q}(T)\geq \rho_{Q}(B(n;n-8,1,1,2)),$
where the equality holds if and only if $T\cong B(n;n-8,1,1,2)$.}

{\bf Proof.} Suppose that $T' \in\mathfrak{T}(n)\cap\Re(n)$ is a graph with minimal distance signless Laplacian spectral radius, where $n\geq8$.
Let $d$ be the diameter of $T'$. Now we verify the following three claims.

{\bf Claim 1.} $d=4$.

Suppose on the contrary that $d\geq5$. Let $P=v_{1}v_{2}\ldots v_{d+1}$ be a path.
Making a C-transformation of $T'$ for the edge $v_{2}v_{3}$ and $v_{d-1}v_{d}$
we obtain the resulting graph $T'_{v_{2}v_{3}}$ and $T'_{v_{d-1}v_{d}}$, respectively.
Obviously, at least one of $T'_{v_{2}v_{3}}$ and $T'_{v_{d-1}v_{d}}$ falls into $\mathfrak{T}(n)\cap\Re(n)$.
By Lemma 2.1, $\rho_{Q}(T')>\rho_{Q}(T'_{v_{2}v_{3}})$ and $\rho_{Q}(T')>\rho_{Q}(T'_{v_{d-1}v_{d}})$.
This contradicts minimality of $T'$, and so $d=4$.

The Claim 1 shows that $T'\cong B(n;n_{0},n_{1},n_{2},\ldots,n_{r})$, where $r\geq 3$.

{\bf Claim 2.} $r=3$.

Suppose on the contrary that $r\geq 4$. Let $T''=B(n;n_{0}+n_{1}+1,n_{2},n_{3},\ldots,n_{r})$.
Thus $T''\in\mathfrak{T}(n)\cap\Re(n)$, and by Lemma 2.1, $\rho_{Q}(T')>\rho_{Q}(T'')$.
This contradiction shows $r=3$, and so $T'\cong B(n;n_{0},n_{1},n_{2},n_{3})$.

{\bf Claim 3.} $T'\cong B(n;n-8,1,1,2)$.

If either one of $n_1$ and $n_2$ is not less than $2$ or $n_{3}\geq 3$, then by Lemma 2.4, $\rho_{Q}(B(n;n_{0},n_{1},$
$n_{2},n_{3}))>\rho_{Q}(B(n;n-8,1,1,2))$, and so $T'\cong B(n;n-8,1,1,2)$. $\Box$

\section{The tree with maximal distance Signless Laplacian spectral radius in $\Re(n)$ or $\mathfrak{T}(n)$}

\qquad{\bf Lemma 3.1.} {\it Let $G=G_{1}\cup G_{2}\cup G_{3}$, where $V(G_{i})\cap V(G_{j})=v_{0}$ for $1\leq i\neq j\leq 3$ and $V(G_{i})\geq 2$ for $i=1,2,3$.
Suppose that $u\in V(G_{2})\backslash \{v_{0}\}$, and
let $\widetilde{G}$ be the graph obtained from $G$ by moving $G_3$ from $v_0$ to $u$.
Let $x=(x_{v_{1}},x_{v_{2}},\cdots ,x_{v_{n}})^T$ be a Perron vector corresponding to $\rho_{Q}(G)$.
If $\sum\limits_{v_{i}\in V(G_{3})\setminus \{v_{0}\}}\sum\limits_{v_{j}\in V(G_{1})}(x_{v_i}+x_{v_j})^{2}\geq \sum\limits_{v_{i}\in V(G_{3})\setminus \{v_{0}\}}\sum\limits_{v_{j}\in V(G_{2})}(x_{v_i}+x_{v_j})^{2}$,
then $\rho_{Q}(\widetilde{G})>\rho_{Q}(G)$.}

{\bf Proof.} It can be easily observed that $d_G(u,v_0)=d_{\widetilde{G}}(u,v_0)$ and that for $v_i\in V(G_{3})\backslash \{v_0\}$, $v_{j}\in V(G_{1})$, $d_{\widetilde{G}}(v_{i},v_{j})=d_{G}(v_{i},v_{j})+d_G(u,v_0)$;
and for $v_{i}\in V(G_{3})\backslash \{v_{0}\}$, $v_{j}\in V(G_{2})$,
$d_{\widetilde{G}}(v_{i},v_{j})\leq d_{G}(v_{i},v_{j})+ d_G(u,v_0)$.

Thus,
$\rho_{Q}(\widetilde{G})-\rho_{Q}(G)\geq x^{T}(Q_{D}(\widetilde{G})-Q_D(G))x$

\qquad\qquad\quad\quad\qquad\quad\,\,\:\,\,$\geq d_G(u,v_0)\sum\limits_{v_{i}\in V(G_{3})\setminus \{v_{0}\}}(\sum\limits_{v_{j}\in V(G_{1})}(x_{v_i}+x_{v_j})^{2}-\sum\limits_{v_{j}\in V(G_{2})}(x_{v_i}+x_{v_j})^{2})$

\qquad\qquad\quad\quad\qquad\quad\,\,\:\,\,$\geq 0$.

Suppose on the contrary that $\rho_{Q}(\widetilde{G})=\rho_{Q}(G)$. Then $x^{T}Q_D(G)x=\rho_{Q}(\widetilde{G})$.

Therefore, $(\rho_{Q}(\widetilde{G})-\rho_{Q}(G))x_{v_0}=\sum\limits_{v_{k}\in V(G_{3})\backslash \{v_{0}\}}d_G(u,v_0)(x_{v_0}+x_{v_k})>0$.
This contradiction shows that $\rho_{Q}(\widetilde{G})>\rho_{Q}(G)$. $\Box$\vskip 2mm

{\bf Lemma 3.2} \cite{H. Liu}. {\it Let $G_{p,q}$ be the connected graph obtained from a graph $G$ of order at least two by attaching two pendant paths of length $p$ and $q$ at a vertex $u$ of $G$.
If $q\geq p\geq 1$, then $\rho_{Q}(G_{p-1,q+1})>\rho_{Q}(G_{p,q})$.} \vskip 2mm

Denote by $T(n,k; t_1,t_2)$ the graph obtained from a path $P_{\ell}$ of order $\ell$ ($\ell\geq 2$) by attaching $t_1\geq 1$ edges at an end of $P_{\ell}$ and $t_2\geq 1$ edges at other end
of $P_{\ell}$, where $t_{1}+t_{2}=k$ and $\ell+k=n$. Clearly, $P_n\cong T(n,2; 1,1)$.
Let $$\mathfrak{D}(n,k)=\{T(n,k; t_1,t_2)|~ t_{1}\geq 1, t_{2}\geq 1~ \mbox{and}~ t_{1}+t_{2}=k\geq 2\}.$$\vskip 2mm

{\bf Theorem 3.3.} {\it The graph with the maximum distance signless Laplacian spectral radius in  $\Re(n)$ lies to $\mathfrak{D}(n,k)$.}

{\bf Proof.} Suppose that $T_0$ is a tree with the maximum distance signless Laplacian spectral radius in $\Re(n)$.
Let $k$ be the number of pendant vertices of $T_{0}$. Then $k\geq 4$.
Let $n^3_{T_{0}}$ be the number of vertices which are of degree at least three in $T_{0}$. Now we verify that the following two claims are true.

{\bf Claim 1.} $n^3_{T_{0}}=2$.

It is clear that $n^3_{T_{0}}\geq 2$. Now suppose by a way of a contradiction that $n^3_{T_{0}}\geq 3$.
Let $A=\{w_{1},w_{2},\ldots,w_{n^3_{T_{0}}}\}$ be the set of all vertices in $T_{0}$ which are of degree at least three,
where $d_{T_{0}}(w_1,w_2)=max\{d_{T_{0}}(w_p,w_q)\mid 1\leq p\neq q\leq n^3_{T_{0}}\}$.
Let $T_{0}=T_{1}\cup T_{2}\cup T_{3}$,
where $V(T_{i})\cap V(T_{j})=w_{3}$ $(1\leq i\neq j\leq3)$ and $w_{1}\in T_{1}$, $w_{2}\in T_{2}$.

If $\sum\limits_{v_{i}\in V(T_{3})\setminus \{w_{3}\}}\sum\limits_{v_{j}\in V(T_{1})}(x_{v_i}+x_{v_j})^{2}\geq \sum\limits_{v_{i}\in V(T_{3})\setminus \{w_{3}\}}\sum\limits_{v_{j}\in V(T_{2})}(x_{v_i}+x_{v_j})^{2}$,
then let $T'_{0}=T_{0}-\sum\limits_{v_{j}\in N_{T_{3}}(w_3)}w_{3}v_{j}+\sum\limits_{v_{j}\in N_{T_{3}}(w_3)}w_{2}v_{j}$,
and otherwise $T'_{0}=T_{0}-\sum\limits_{v_{j}\in N_{T_{3}}(w_3)}w_{3}v_{j}+\sum\limits_{v_{j}\in N_{T_{3}}(w_3)}w_{1}v_{j}$.
Thus $T'_{0}\in\Re(n)$, and by Lemma 3.1, $\rho_{Q}(T'_{0})>\rho_{Q}(T_{0})$.
This contradicts maximality of $T_{0}$, and so $n^3_{T_{0}}=2$.\vskip 2mm

The Claim 1 shows that $T_{0}$ consists of a path with its two end vertices $w_1$ and $w_2$ attached some pendant paths respectively.\vskip 2mm

{\bf Claim 2.} $T_{0}\cong T(n,k; t_1,t_2)$.

Suppose on the contrary that $T_{2}=w_{1}u_{1} \ldots u_{s}$ is a pendant path of length $s\geq 2$.
Let $T_{0}=T_{1}\cup T_{2}\cup T_{3}$ be shown in Fig. 3,
where $V(T_{i})\cap V(T_{j})=w_1$ for $1\leq i\neq j\leq 3$, $N_{T_0}(w_{1})\backslash\{u_{1},v_{1}\} \subset V(T_{1})$ and $w_{2}\in V(T_3)$.
\begin{center}
\begin{figure}[htbp]
\centering
\includegraphics[height=20mm]{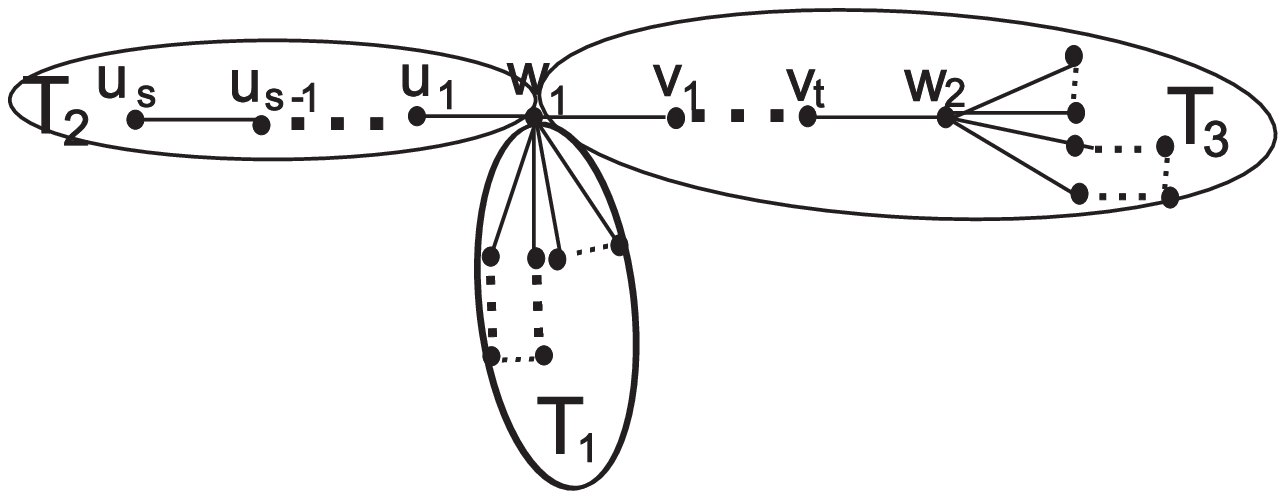}
\\ Fig. 3. ~~~~ $T_0$
\end{figure}
\end{center}

If $\sum\limits_{v_{i}\in V(T_{1})}\sum\limits_{v_{j}\in V(T_{3})}(x_{v_i}+x_{v_j})^{2}\geq \sum\limits_{v_{i}\in V(T_{1})}\sum\limits_{v_{j}\in V(T_{2})}(x_{v_i}+x_{v_j})^{2}$,
then let $\widetilde{T_{0}}=T_{0}-\sum\limits_{v_{j}\in V(T_{1})}w_{1}v_{j}+\sum\limits_{v_{j}\in V(T_1)}u_{s-1}v_j$,
and otherwise $\widetilde{T_{0}}=T_{0}-\sum\limits_{v_{j}\in V(T_{1})}w_{1}v_{j}+\sum\limits_{v_{j}\in V(T_1)}w_{2}v_j$.
By Lemma 3.1, $\rho_{Q}(\widetilde{T_{0}})>\rho_{Q}(T_{0})$.
This contradicts maximality of $T_0$, and so $T_{0}\cong T(n,k; t_1,t_2)$.

Let $n_{1},n_{2},\ldots,n_{r}$ be positive integers satisfying $n_{1}\leq n_{2}\leq \ldots\leq n_{r}$ and $\sum\limits_{i=1}^{r}n_{i}=n-1$.
Then we denote by $S(n;n_{1},n_{2},\ldots,n_{r})$ the tree obtained by adding an edge between a fixed vertex $u$ and an end vertex of the path $P_{n_i}$ for each $1\leq i\leq r$. \vskip 2mm

{\bf Theorem 3.4.} {\it If $T \in\mathfrak{T}(n)$, then $\rho_{Q}(S(n;2,2,n-5))\geq\rho_{Q}(T)$, where the equality holds
if and only if $T\cong S(n;2,2,n-5)$.}

{\bf Proof.}  It is clear that there is such a vertex $v$ in $T$ that $N(v)$ contains at least three non-pendant vertices.
Thus, if $n_1, n_2,\ldots,n_r$ are respectively orders of $r$ connected components of $T-v$ then at least three of them are not less than two.
By Lemma 3.2, $\rho_{Q}(S(n;n_{1},n_{2},\ldots,n_{r}))\geq \rho_{Q}(T)$.
We can further see by Lemma 3.2 that $\rho_{Q}(S(n;2,2,n-5))\geq \rho_{Q}(S(n;n_{1},n_{2},\ldots,n_{r})).$   $\Box$\vskip 2mm

Let~$\mathfrak{B}(n,k)=\{T|~T~ \mbox{is a tree of order}~n~ \mbox{with}~ k ~\mbox{pendant vertices}, \mbox{where}~2 \leq k \leq n-2\}$.
Let~$P(n;i,j)$ be the tree obtained from the path $P=v_{1}v_{2}\ldots v_{n-3}$ by attaching a pendant edge at $v_{i}$ and a pendant
path of length $2$ at $v_{j}$.\vskip 2mm

{\bf Lemma 3.5.} {\it Let $T \in\mathfrak{T}(n)\cap\Re(n)\cap \mathfrak{B}(n,4)$. Then we have $$max\{\rho_{Q}(P(n;2,3)), \rho_{Q}(P(n;2,n-5))\}\geq \rho_{Q}(T),$$
where the equality holds if and only if $T\cong P(n;2,3)$ or $T\cong P(n;2,n-5)$.}

{\bf Proof.} Suppose that $T_{0} \in\mathfrak{T}(n)\cap\Re(n)\cap \mathfrak{B}(n,4)$ is a tree with the maximum distance signless Laplacian spectral radius.
Let $d$ be the diameter of $T_{0}$, and $P=u_{1}u_{2}\ldots u_{d+1}$ be a path in $T_{0}$.
Then there exist one pendant path $P^{1}$ of length $p\geq 1$ at some $u_{i}$ and another one pendant path $P^{2}$ of length $q\geq 2$ at another vertex $u_{j} (\neq u_i)$.
We assume without loss of generality that $q\geq p$. Now we verify the following three claims.

{\bf Claim 1.} $p=1$ and $q=2$.

Set $P^{1}=u_{i}z_{1}z_{2}\ldots z_{p}$ and $P^{2}=u_{j}w_{1}w_{2}\ldots w_{q}$.
Suppose on the contrary that $p\geq 2$.
Note that $p\leq i-1$. Hence we let $T_{0}'=T_{0}-z_{p-1}z_{p}+u_{1}z_{p}$.
Clearly, $T_{0}' \in\mathfrak{T}(n)\cap\Re(n)\cap \mathfrak{B}(n,4)$, and by Lemma 3.2, $\rho_{Q}(T_{0}')>\rho_{Q}(T_{0})$.
This contradiction shows that $p=1$.

Similarly, we can prove that $q=2$.

The Claim 1 shows that $T_{0}\cong P(n;i,j)$, where $2\!\leq\!i\!<\!j\!\leq\!n-5$ and $P=u_{1}u_{2}\ldots u_{n-3}$.

{\bf Claim 2.} $i=2$.

Suppose that $u_{n-2}$ and $u_{n-1}$ are respectively the neighbors of $u_{i}$ and $u_{j}$ in $T_{0}$.
Let $x$ be a Perron vector corresponding to $\rho_{Q}(T_{0})$.
Suppose on the contrary that $i\geq 3$. Let $U_{1}=\{u_{1},u_{2},\ldots,u_{i-1}\}$, $U_{2}=\{u_{i+1},u_{i+2},\ldots,u_{n-3},u_{n-1},u_{n}\}$ and $U_{3}=\{u_{i},u_{n-2}\}$.
If $\sum\limits_{v_{i}\in U_{3}\setminus \{u_{i}\}}\sum\limits_{v_{j}\in U_{1}}(x_{v_i}+x_{v_j})^{2}\geq \sum\limits_{v_{i}\in U_{3}\setminus \{u_{i}\}}\sum\limits_{v_{j}\in U_{2}}(x_{v_i}+x_{v_j})^{2}$,
then let $T_{0}'=T_{0}-u_{i}u_{n-2}+u_{n-1}u_{n-2}$, and otherwise $T_{0}'=T_{0}-u_{i}u_{n-2}+u_{i-1}u_{n-2}$.
Thus $T_{0}' \in\mathfrak{T}(n)\cap\Re(n)\cap \mathfrak{B}(n,4)$, and by Lemma 3.1, $\rho_{Q}(T_{0}')>\rho_{Q}(T_{0})$.
This contradicts maximality of $T_{0}$, and so $i=2$.

The Claim 2 shows that $T_{0}\cong P(n;2,j)$.

{\bf Claim 3.} $T_{0}\cong P(n;2,3)$ or  $T_{0}\cong P(n;2,n-5)$.

Suppose on the contrary that $T_{0}\ncong P(n;2,3)$ and $T_{0}\ncong P(n;2,n-5)$. Then $4\leq j\leq n-6$.
Let $U_{1}=\{u_{1},u_{2},\ldots,u_{j-1},u_{n-2}\}$, $U_{2}=\{u_{j+1},u_{j+2},\ldots,u_{n-3}\}$ and $U_{3}=\{u_{j},u_{n-1},u_{n}\}$.
If $\sum\limits_{v_{i}\in U_{3}\setminus \{u_{j}\}}\sum\limits_{v_{j}\in U_{1}}(x_{v_i}+x_{v_j})^{2}\geq \sum\limits_{v_{i}\in U_{3})\setminus \{u_{j}\}}\sum\limits_{v_{j}\in U_{2}}(x_{v_i}+x_{v_j})^{2}$,
then let $T_{0}''=T_{0}-u_{j}u_{n-1}+u_{j+1}u_{n-1}$, and otherwise $T_{0}''=T_{0}-u_{j}u_{n-1}+u_{j-1}u_{n-1}$.
Thus $T_{0}'' \in\mathfrak{T}(n)\cap\Re(n)\cap \mathfrak{B}(n,4)$, and by Lemma 3.1, $\rho_{Q}(T_{0}'')>\rho_{Q}(T_{0})$.
This contradicts maximality of $T_{0}$, and so $T_{0}\cong P(n;2,3)$ or $T_{0}\cong P(n;2,n-5)$. $\Box$\vskip 2mm

Suppose that $P=v_1v_2\ldots v_s$ is a path of a graph $G$ and $u\in V(G)$. Then we let $d_G(u,P)=min\{d_G(u,v_i)|1\leq i\leq s\}.$ \vskip 2mm

{\bf Theorem 3.6.} {\it Let $T \in\mathfrak{T}(n)\cap\Re(n)$. Then we have $$max\{\rho_{Q}(P(n;2,3)), \rho_{Q}(P(n;2,n-5))\}\geq \rho_{Q}(T),$$
where the equality holds if and only if $T\cong P(n;2,3)$ or $T\cong P(n;2,n-5)$.}

{\bf Proof.} Suppose that $T_{0} \in\mathfrak{T}(n)\cap\Re(n)$ is a tree with the maximum distance signless Laplacian spectral radius.
Let $A=\{w_{1},w_{2},\ldots,w_{n^3_{T_{0}}}\}$ be the set of all vertices in $T_{0}$ which are of degree at least three.
It is obvious that $n^3_{T_{0}}\geq 2$.
Let $P=u_{1}u_{2}\ldots u_{d+1}$ be a path, where $d$ is the diameter of $T_{0}$. Next we distinguish two cases.

{\bf Case 1.} $|A\cap V(P)|\geq 2$.

Suppose for a contradiction that $A\nsubseteq V(P)$.
Without loss of generality we assume that $w_{\ell}\in A\setminus V(P)$ such that $d_{T_{0}}(w_{\ell},P)=max\{d_{T_{0}}(w_i,P)\mid 1\leq i\leq n^3_{T_{0}}\}$.
Then there are at least two pendant paths at $w_{\ell}$, say, $P^{1}=w_{\ell}z_{1} \ldots z_{p}$ and $P^{2}=w_{\ell}v_{1} \ldots v_{q}$ ($p\geq q\geq1$).
Let $T_{0}'=T_{0}-v_{q-1}v_{q}+z_{p}v_{q}$. Since $|A\cap V(P)|\geq 2$, we have $|A\backslash\{w_{\ell}\}|\geq 2$, and so $T_{0}'\in\mathfrak{T}(n)\cap\Re(n)$.
By Lemma 3.2, $\rho_{Q}(T_{0}')>\rho_{Q}(T_{0})$.
This contradicts maximality of $T_{0}$, and so $A\subseteq V(P)$.

Suppose by a way of contradiction that $n_{T_{0}}^{3}\geq 3$.
Let $d_{T_{0}}(w_1,w_2)=max\{d_{T_{0}}(w_p,w_q)\mid 1\leq p\neq q\leq n^3_{T_{0}}\}$.
Then there are two pendant paths at $w_{1}$, say, $P^{3}=w_{1}y_{1} \ldots y_{s}$ and $P^{4}=w_{1}f_{1} \ldots f_{t}$ ($s\geq t\geq1$).
Let $T_{0}'=T_{0}-f_{t-1}f_{t}+y_{s}f_{t}$.
Thus $T_{0}'\in\mathfrak{T}(n)\cap\Re(n)$,
and by Lemma 3.2, $\rho_{Q}(T_{0}')>\rho_{Q}(T_{0})$.
This contradicts maximality of $T_{0}$, and so $n^3_{T_{0}}=2$.

Suppose for a contradiction that $d_{T_{0}}(w_{1})\geq 4$.
Then there exist three pendant paths at $w_{1}$ in $T_{0}$, say $P^{3}=w_{1}y_{1} \ldots y_{s}$, $P^{4}=w_{1}f_{1} \ldots f_{t}$ and $P^{5}=w_{1}h_{1} \ldots h_{m}$ ($s\geq t\geq m\geq 1$).
Let $T_{0}'=T_{0}-h_{m-1}h_{m}+f_{t}h_{m}$.
Thus $T_{0}'\in\mathfrak{T}(n)\cap\Re(n)$, and by Lemma 3.2, $\rho_{Q}(T_{0}')>\rho_{Q}(T_{0})$.
This contradicts maximality of $T_{0}$, and so $d_{T_{0}}(w_{1})=3$.
Similarly, we can prove $d_{T_{0}}(w_{2})=3$.
Therefore, $T_{0} \in\mathfrak{T}(n)\cap\Re(n)\cap \mathfrak{B}(n,4)$, and
from Lemma 3.5 we know that the assertion is true.

{\bf Case 2.} $|A\cap V(P)|=1$.

Let $A\cap V(P)=\{w_{1}\}$.
Since $n^3_{T_{0}}\geq 2$, $A\backslash V(P)\neq \varnothing$.
Suppose on the contrary that $n^3_{T_{0}}\geq 3$.
Let $w_{2}\in A\backslash \{w_{1}\}$ such that $d_{T_{0}}(w_1,w_2)=max\{d_{T_{0}}(w_1,w_j)\mid 2\leq j\leq n^3_{T_{0}}\}$.
Then there are two pendant paths at $w_{2}$, say, $P^{6}=w_{2}a_{1} \ldots a_{r}$ and $P^{7}=w_{2}b_{1} \ldots b_{j}$ ($r\geq j\geq1$).
Let $T_{0}''=T_{0}-b_{j-1}b_{j}+a_{r}b_{j}$.
Since $|A\backslash\{w_{2}\}|\geq 2$, $T_{0}''\in\mathfrak{T}(n)\cap\Re(n)$, and by Lemma 3.2, $\rho_{Q}(T_{0}'')>\rho_{Q}(T_{0})$.
This contradicts maximality of $T_{0}$, and so $n^3_{T_{0}}=2$.

By similar argument as Case 1, we can verify that $d_{T_{0}}(w_{1})=d_{T_{0}}(w_{2})=3$,
and so $T \in\mathfrak{T}(n)\cap\Re(n)\cap \mathfrak{B}(n,4)$.
By Lemma 3.5, we know that the assertion is true.     $\Box$

\end{document}